\documentclass[12pt,journal]{IEEEtran}
\pdfoutput=1
\usepackage{cite}

\ifCLASSINFOpdf
   \usepackage[pdftex]{graphicx}
\else
\fi
\usepackage{amssymb}
\usepackage{url}

\begin{document}
\title{Numerical Reproducibility and Parallel Computations: Issues for Interval Algorithms}

\author{Nathalie~Revol,~\IEEEmembership{Member,~IEEE,}
        and Philippe Th\'eveny%
\thanks{Universit\'e de Lyon - LIP (UMR 5668 CNRS - ENS de Lyon - INRIA - UCBL), ENS de Lyon, France}%
\thanks{Nathalie Revol is with
INRIA -
\texttt{Nathalie.Revol@ens-lyon.fr} - WWW home page:
\texttt{http://perso.ens-lyon.fr/nathalie.revol/}}%
\thanks{Philippe Th\'eveny is with
ENS de Lyon -
\texttt{Philippe.Theveny@ens-lyon.fr} - WWW home page:
\texttt{http://perso.ens-lyon.fr/philippe.theveny/}
}
}

\markboth{Numerical Reproducibility and Interval Algorithms}{author name to be inserted}

\maketitle

%
\IEEEpeerreviewmaketitle

\begin{abstract}
What is called \emph{numerical reproducibility} is the problem of getting the same
result when the scientific computation is run several times, either on the
same machine
or on different machines,
with different types and numbers of processing units, execution environments,
computational loads etc.
This problem is especially stringent for HPC numerical simulations.
In what follows, the focus is on parallel implementations of interval arithmetic using floating-point arithmetic.
For interval computations, numerical reproducibility is of course
an issue for testing and debugging purposes. However, as long as the computed
result encloses the exact and unknown result, the inclusion property, which is
the main property of interval arithmetic, is satisfied and getting bit for bit
identical results may not be crucial.
Still, implementation issues may invalidate the inclusion property.
Several ways to preserve the inclusion property are presented, on the example of the
product of matrices with interval coefficients.
\end{abstract}

\begin{IEEEkeywords}
interval arithmetic, numerical reproducibility, parallel implementation, floating-point arithmetic, rounding mode.
\end{IEEEkeywords}

\section{Introduction}
Interval arithmetic is a mean to perform numerical computations and to get a guarantee
on the computed result. In interval arithmetic, one computes with intervals, not
numbers. These intervals are guaranteed to contain the (set of) exact values,
both input values and computed values, and this property is preserved throughout
the computations. Indeed, the fundamental theorem of interval arithmetic is
the \emph{inclusion property}: each computed interval contains the exact (set of)
result(s).
To learn more about interval arithmetic, see \cite{Moore66:IA,MoKeCl09:IA,Neumaier90:IA,Rump05:IA,Rump10:IA,Tucker11:IA}.

The inclusion property is satisfied by the definition of the arithmetic operations,
and other operations and functions acting on intervals. Implementations of interval
arithmetic on computers often rely on the floating-point unit provided by the
processor. The inclusion property is preserved via the use of directed roundings.
If an interval is given exactly for instance by its endpoints, the left endpoint
of its floating-point representation
is the rounding towards $-\infty$ of the left endpoint of the exact interval.
Likewise, the right endpoint of its floating-point representation is the rounding towards $+\infty$ of the right endpoint 
of the exact interval. Similarly, directed roundings are used to implement mathematical
operations and functions: in this way, roundoff errors are accounted for and the inclusion property is satisfied.

In this paper, the focus is on the problems encountered while implementing interval
arithmetic, using floating-point arithmetic, on multicore architectures. 
In a sense these issues relate to issues known as problems of
\emph{numerical reproducibility} in scientific computing using floating-point arithmetic
on emerging architectures. The common points and the differences and concerns which are specific
to interval arithmetic will be detailed.

The main contributions of this paper are \textbf{the identification of problems that cause numerical irreproducibility
of interval computations,
and recommendations to circumvent these problems}.
The classification proposed here distinguishes between three categories. 
The first category, addressed in Section \ref{Sec.Prec}, concerns problems of variable computing precision
that occur both in sequential and parallel implementations, both for floating-point and interval computations.
These behaviours are motivated by the quest of speed, at the expense of accuracy on the results, be they
floating-point or intervals results.
However, even if these behaviours hinder numerical reproducibility, usually they do not threaten
the validity of interval computations: interval results satisfy the inclusion property.

The second category, developed in Section \ref{Sec.Order}, concerns the order of operations.
Due mainly to the indeterminism in multithreaded computations, the order of operations may vary.
It is the most acknowledged problem.
Again, this problem is due to the quest of short execution time, again it hinders numerical reproducibility
of floating-point or interval computations, but at first sight it does not invalidate the inclusion property.
However, an interval algorithm will be presented, whose validity relies on an assumption on the order
of operations.

The problems of the third category, detailed in Section \ref{Sec.Rndg}, are specific to interval computations
and have an impact on the validity of interval results. The question is whether directed rounding modes,
set by the user or the interval library, are respected by the environment.
These problems occur both for sequential and parallel computations. 
Ignoring the rounding modes can permit to reduce the execution time.
Indeed, changing the rounding mode can incur a severe penalty: the slowdown lies between 10 and 100
on architectures where the rounding modes are accessed via global flags and where changing a rounding mode
implies to flush pipelines. However, it seems that the most frequent motivation is either
the ignorance of issues related to rounding modes, or the quest for an easier and faster development of the execution environment
by overlooking these issues.

For each category, 
we give recommendations for the design of interval algorithms:
numerical reproducibility is still not attained, but the impact of the problems identified
is attenuated and the inclusion property is preserved. 
\\

Before detailing our classification, numerical reproducibility is defined in Section \ref{Sec.Pb} and
then an extensive bibliography is detailed in Section \ref{Sec.Prev}.

\section{Problem of Numerical Reproducibility}
\label{Sec.Pb}
In computing in general and in scientific computing in particular, the quest for speed
has led to parallel or distributed computations, from multithreaded programming to high-performance computing (HPC).
Within such computations, contrary to sequential ones, the order in which events occur is not deterministic.
Events here is a generic term which covers communications, threads or processes, modification of 
shared data structures\ldots.
The order in which communications from different senders arrive at a common receiver,
the relative order in which threads or processes are executed, the order in which shared data structures
such as working lists are consulted or modified, is not deterministic.
Indeed, for the example of threads, the quest for speed leads a greedy scheduler to start a task
as soon as it is ready, and not at the time or order at which it should start in a sequential execution.
One consequence of this indeterminism is that speed is traded against numerical reproducibility of
floating-point computations: computed results may differ, depending on the order in which events
happened in one run or another.
However, numerical reproducibility is a desirable feature, at least for debugging purposes:
how is it possible to find and fix a bug that occurs sporadically, in an indeterministic way?
Testing is also impossible without numerical reproducibility: when a program returns results
which differ from the results of a reference implementation, nothing can be concluded
in the absence of numerical reproducibility.
To address the issue of numerical reproducibility,
in \cite{StBaBoLVRiSt13:repro}, not only the use of higher computing precision is recommended,
but also the use of interval arithmetic or of some variant of it in order to get numerical guarantees,
and the use of tools that can diagnose sensitive portions of code and take corrective actions.
We follow these tracks and we focus on interval arithmetic, the same conclusions can be drawn for its 
``variants'' such as affine arithmetic or polynomial models, e.g. Taylor models.

First, let us introduce the problem of numerical reproducibility in floating-point arithmetic.
Floating-point arithmetic is standardised and the IEEE-754 standard \cite{IEEE-754-85,IEEE-754-08}
specifies completely the formats of floating-point numbers, their binary representation and the
behaviour of arithmetic operations. One of the goals of this standardisation was to enable portability
and reproducibility of numerical computations across architectures.
However, several runs of the same code on different platforms, or even on the same platform when the code
is parallel (typically, multithreaded) may yield different results.
Getting the same result on the same platform is sometimes called \emph{repeatability}. Only the word
\emph{reproducibility} will be used throughout this paper.

Let us give a first explanation to the lack of numerical reproducibility.
On some architectures,
and with some programming languages and compilers, intermediate computations such as the intermediate
sum in $s \leftarrow a+b+c$ can take place with variable precision. This precision can depend on
the availability of extended size registers (such as 80-bits registers on x87 or IA64 architectures,
when computing with 64-bit operands),
on whether operations use registers or store intermediate results in memory,
on data alignment or on the execution of other threads \cite{ChGoRaAhLe13:repro-HPC},
on whether an intermediate computation is promoted to higher precision or not. 
The impact of the computing precision on numerical
irreproducibility is discussed in details in Section \ref{Sec.Prec}.

Let us stick to our toy example $s \leftarrow a+b+c$ to illustrate a second
cause to numerical irreproducibility. As floating-point addition is not associative
(neither is multiplication), the order according to which the operations are performed matters.
Let us illustrate this with an example in double precision: $\mbox{RN} (1+(2^{100} - 2^{100}))$ (where $\mbox{RN}$ stands
for \emph{rounding-to-nearest}) yields $1$, since $\mbox{RN}(2^{100}-2^{100}) = 0$,
whereas $\mbox{RN} ((1+2^{100})- 2^{100})$ yields $0$, since $\mbox{RN}(1+2^{100}) = 2^{100}$.
The last result is called a \emph{catastrophic cancellation}: most or all accuracy is lost when close (and even equal in this example)
large numbers are subtracted and only roundoff errors remain, hiding the meaningful result.
Real-life examples of this problem are to be found in \cite{HeDi01:repro}, where the results of computations
for ocean-atmosphere simulation strongly depends on the order of operations.
Other examples are given in \cite{Diethelm12:repro}
for the numerical simulation of punching of metal sheets.
This lack of associativity of floating-point operations implies that the result of a reduction
operation depends on the order according to which the operations are performed.
However, with multithreaded or parallel or distributed implementations, this order is not deterministic and 
the computed result thus varies from one execution to the next.
How the order of operations influences numerical computations and interval computations is detailed
in Section \ref{Sec.Order}.

A last issue, which is specific to the implementation of interval arithmetic on parallel environments,
is the respect
of rounding modes. As already mentioned, the implementation of interval arithmetic crucially requires
directed rounding modes.
However, it has been observed that rounding modes are not respected by numerical libraries such
as BLAS \cite{LaMeMo12}, nor by compilers when default options are used,
 nor by execution environment for multithreaded computations, nor by parallel
languages such as OpenMP or OpenCL. Either this is not documented, or this is explicitly mentioned as
being not supported, cf. Section \ref{Sec.Rndg} for a more detailed discussion.
Respecting the rounding modes is required by the IEEE-754 standard for floating-point arithmetic
and the behaviours just mentioned are either misbehaviours as in the example of Section \ref{SubSec.Rndg.Pbs},
or are (often undocumented) features of the libraries, often "justified" by the quest of shorter execution times.
\\

These phenomena explain the lack of numerical reproducibility for floating-point computations,
i.e. the fact that two different runs yield two different results, or the loss of the inclusion
property in interval arithmetic, due to the non-respect of rounding modes.
\\

Facing this lack of reproducibility, various reactions are possible.
One consists in acknowledging the computation of differing results as
an indication of a lack of numerical stability. In a sense, a positive way
of considering numerical irreproducibility is to consider it as useful information on the numerical
quality of the code.
However, for debugging and testing purposes, reproducibility is more than helpful.
For such purposes, \emph{numerical reproducibility means getting bitwise identical results from one run to the next}.
Indeed, this is the most common definition of numerical reproducibility. This definition is also
useful for contractual or legal purposes (architectural design, drug design are instances mentioned in the slides corresponding to \cite{DeNg13:repro}),
as long as a reference implementation, or at least a reference result, is given.
However, this definition is not totally satisfactory as the result is not well-defined.
In particular, it says nothing about the accuracy of the result.
Requiring the computed result to be the correct rounding of the exact result
is a  semantically meaningful definition of  numerical reproducibility.
The computed result is thus uniquely defined.
The difficulty with this definition is to devise an algorithm that is efficient on parallel platforms
and that computes the correct rounding of the exact result.
This definition has been adopted in \cite{MIScDi06:LHC}, for LHC computations of 600,000 jobs
on 60,000 machines: efficient mathematical functions such as exponential and logarithm, with correct rounding, were available.
However, in most cases it is not known how to compute efficiently the correctly rounded result.
For this reason, this definition of numerical reproducibility as computing the correct rounding of the exact result
may be too demanding. Thus we prefer to keep separate the notions of numerical reproducibility (i.e. getting
bitwise identical results) and of correct rounding.
\\

To sum up, numerical reproducibility can be defined as getting bitwise identical results, where these results
have to been specified in some way, e.g. by a reference implementation which can be slow.
We consider it is simpler not to mix this notion with the notion of correct rounding of the exact result,
which is uniquely specified.
Our opinion has been reinforced
by the results in \cite{DeNg13:repro}: it is even possible to define several levels of numerical reproducibility,
each one corresponding to a level of accuracy. One has thus a hierarchy of reproducibilities, corresponding to
different tradeoffs between efficiency and accuracy, ranging from low numerical quality to correct rounding.

\section{Previous Work}
\label{Sec.Prev}
In the previous section, we introduced various sources of numerical irreproducibility 
and we delineated their main features. More detailed and technical explanations are given
in \cite{Gropp05}: 
implementation issues such as data race, out-of-order execution, message buffering with
insufficient buffer size, 
non-blocking communication operations are also introduced.
A tentative classification of indeterminism can be found in \cite{ChGoRaAhLe13:repro-HPC},
it distinguishes between \emph{external determinism} that roughly corresponds to getting the
same results independently of the internal states reached during the computation,
and \emph{internal determinism} that requires that internal execution steps are the same
from one run to the other.
External indeterminism can be due to data alignment that varies from one execution to the next,
order of communications with ``wildcard receives" of MPI, and other causes already mentioned.
In what follows, only external indeterminism  will be considered.
\\

Numerical irreproducibility, in particular of summations of floating-point numbers,
is cited as early as 1994 in \cite{KaOiSaSeHi94} for
weather forecasting applications. More recently, it has been clearly put in evidence in \cite{Diethelm12:repro}
where the application is the numerical simulation of a deep drawing process for deforming metal sheets:
depending on the execution, the computed variations of the thickness of the metal sheet vary.
Other references mentioning numerical irreproducibility are for instance \cite{HeDi01:repro}
for an application in ocean-atmosphere simulation, \cite{LeRaWaYa12}
and reference $[10]$ herein about digital breast tomosynthesis,
or \cite{ViCMGuMaKr09:repro} for power state estimation in the electricity grid.
\\

\subsection{The Example of the Summation} 
The non-associativity of floating-point operations is the major explanation to the phenomenon of
numerical irreproducibility. The simplest problem that exemplifies the non-associativity
is the summation of $n$ floating-point numbers. It is also called a reduction of $n$ numbers
with the addition operation. Not surprisingly, many efforts to
counteract numerical irreproducibility focus on the summation problem,
as accuracy can be lost due to catastrophic cancellation or catastrophic absorption\footnote{Catastrophic absorption
occurs in a sum when small numbers are added to a large number and ``disappear'' in the roundoff error, i.e. they do
not appear in the final result, whereas adding first these small numbers would result in a number large enough to be added
to the first large number and be ``visible''.} in summation. Let us list some
of them in chronological order of publication. An early work \cite{HeDi01:repro} on the summation uses
and compares several techniques to get more accuracy on the result: the conclusion is that
compensated summation and the use of double-double arithmetic give the best results.
Following this work, for physical simulations,
in \cite{RoRoAu11:repro} conservation laws were numerically enforced using compensated sums, either
Kahan's version or Knuth's version, and implemented as MPI reduction operators.
Similarly, in \cite{TaPaSaPa10-repro}, a simplified form of ``single-single'' arithmetic
is employed on GPU to sum the energy and satisfy numerically the conservation law.
Bailey, the main author of QD, a library for double-double and quad-double arithmetic \cite{HiLiBa:QD},
also advocates the use of this higher-precision arithmetic in \cite{Bai12:repro}.
However, even if the use of extra computing precision yields accurate results
for more ill-conditioned summations, which means the obtention of results with more
correct bits, it does not ensure numerical reproducibility, as solving even worse-conditioned
problems shows.
Since the problem can be attributed to the variations in the summation order, a first
approach to get reproducible results consists in fixing the
reduction tree \cite{ViCMGuMaKr09:repro}. More precisely, using a fixed integer $K$,
the array of summands is split into $K$ chunks of consecutive
subarrays,
each subarray (or chunk) being summed sequentially (but each chunk can be summed independently of the other
ones) and then the order of the reduction of the $K$ partial sums is also sequential.

Another approach in \cite{DeNg13:repro,DeNg13:repro-HPC} consists in what is called
\emph{pre-rounding}. Even if the actual implementation is very different from the process
explained here, it can be thought of as a sum in fixed-point, as learnt in elementary school.
The mantissa of every summand is aligned with respect to the point, then a leftmost ``vertical
slice'' is considered and added to produce a sum $S_1$. The width of this slice is chosen
in such a way that the sum $S_1$ can be performed exactly using the width of the mantissa
of floating-point numbers (e.g. 53 bits for the double precision format).
As this sum $S_1$ is exact, it is independent of the order in which the additions are
performed and is thus numerically reproducible on any platform.
To get a more accurate result, a second vertical slice, just right to the first one,
can be summed, again yielding an
exact sum $S_2$ and the final result is the (floating-point, thus inexact) sum of
$S_1$ and $S_2$.
An increase of the number $K$ of slices corresponds to an increase of the computing precision.
However, for a fixed
$K$, each partial sum $S_1, \ldots , S_K$ is exact and thus numerically reproducible.
The details to determine the slices and to get $S_1, \ldots S_K$ are given in \cite{DeNg13:repro}
and a faster algorithm for exascale computing, i.e. for really large platforms, is provided in \cite{DeNg13:repro-HPC}.

\subsection{Approaches to Reach Numerical Reproducibility}
In \cite{KrBiMuRe04}
a tool, called MARMOT, that detects and signals race conditions and deadlocks in MPI
codes is introduced, but reduction operations are not handled.

The original proposal for Java by Gosling \cite{Sun01}
included numerical reproducibility of floating-point computations. To reach reproducibility,
it prohibited the use of any format different from Binary32 and Binary64, the use of the FMA
as well as any optimisation based on the associativity of the operators.
It also forbade changes of the rounding modes
\cite{Philippsen98} and the only rounding mode is to-nearest \cite[Section 4.2.4]{Java13}.
It did not include the handling of exceptions via flags, as required by the IEEE-754 standard.
The seminal talk by Kahan in 1998 \cite{KaDa98:Java} has shaken these principles.
Actually Kahan disputed mainly the lack of exception handling and the restriction to shorter (thus, less precise) formats
even when longer ones are available.
It seems that this dispute opened the door to variations around Java. Indeed Java Grande
\cite{PhBoGePoMoGaFo01,Thiruvathukal02}
proposes to allow the use of longer formats and of FMA. The use of the associativity of operations to optimise
the execution time has been under close scrutiny for a longer lapse of time and remains prohibited \cite[Section 15.7.3]{Java13},
as explicit rules are given, e.g. left-to-right priority for $+$, $-$, $\times$.
However, strict adherence to the initial principles of Java can be enforced by using the \texttt{StrictFp} keyword.
For instance, the VSEit environment for modelling complex systems, simulating them and getting a graphical view
\cite{Brassel01:repro}, uses the \texttt{StrictMath} option in Java 
as a way to get numerical reproducibility. However, the need for  getting reproducibility is not explained in much details.

Finally, 
Intel MKL (Math Kernel Library) 11.0 introduces a feature called Conditional Numerical Reproducibility (CNR)
\cite{Intel-CNR12}
which provides functions for obtaining reproducible floating-point results. When using these new features, Intel MKL functions are designed to return the same floating-point results from run-to-run, subject to the following limitations:
\begin{itemize}
\itemsep -1pt
\item calls to Intel MKL occur in a single executable
\item input and output arrays in function calls must be aligned on 16, 32, or 64 byte boundaries on systems with SSE / AVX1 / AVX2 instructions support (resp.)
\item the number of computational threads used by the library remains constant throughout the run.
\end{itemize}
These conditions are rather stringent.
Another approach to numerical reproducibility consists in providing correctly rounded functions, at least for the mathematical library 
\cite{MIScDi06:LHC}.

The approaches presented here are not yet entirely satisfactory, either because they do not really offer
numerical reproducibility or because they handle only summation, or because performances are too drastically
slowed down. Furthermore, none addresses interval computations.
In what follows, we propose a classification of the sources of numerical irreproducibility for interval computations
and some recommendations to circumvent these problems, even if we do not have the definitive solution.
In our classification, a first source of problem is the variability of the employed computing precision.

\section{Computing Precision}
\label{Sec.Prec}

\subsection{Problems for Floating-Point Computations}
\label{SubSec.Prec.PbsFP}
The computing precision used for floating-point computations depends on the employed format
(as defined by IEEE-754 standard \cite{IEEE-754-85,IEEE-754-08}).
The single precision corresponds to the representation of floating-point numbers on 32 bits,
where 1 bit is used for the sign, 8 bits for the exponent and the rest for the significand.
The corresponding format is called \texttt{Binary32}. The double precision uses 64 bits,
hence the name \texttt{Binary64}, with 1 bit for the sign and 11 for the exponent.
The quadruple precision, also known as \texttt{Binary128}, uses 128 bits, with 1 bit for the sign
and 15 bits for the exponent.
The rest of Section \ref{SubSec.Prec.PbsFP} owes much to \cite{dD13:repro}.

On some architectures (IA32 / x87), registers have a longer format: 80 bits instead of 64.
The idea prevailing to the introduction of these long registers
was to provide higher accuracy by computing, for a while, with higher precision than the Binary64 format and to round
the result into a Binary64 number only after several operations. However it entailed the so-called
``double-rounding'' problem: an intermediate result is first rounded into a 80-bit number, then into
a 64-bit number when it is stored, but these two successive roundings can yield a result different
from rounding directly into a 64-bit number.
From the point of view of numerical reproducibility, the use of extended precision registers
is also troublesome, as the operations which take place within registers and the temporary storage
into memory can
occur at different stages of the computation,
they may vary from run to run, depending for instance on the load of the
current processor,
on data alignment or on the execution of other threads \cite{ChGoRaAhLe13:repro-HPC}.

Another issue is the format chosen for the intermediate result, say for the intermediate sums
in the expression $a+b+c+d$, where $a$, $b$, $c$ and $d$ are Binary32 floating-point formats.
Notwithstanding the order in which the intermediate sums are computed, 
let us focus on the problem of the precision.
If the architecture offers Binary64 and if the language is C or Python, then the intermediate
sums may be promoted to Binary64. It will not be the case if the language is Java with the
\texttt{StrictFp} keyword, or Fortran.
In C, it may or may not be the case, depending on the compiler and on the compilation options:
the compiler may prefer to take advantage of a vectorised architecture like SSE2 or AVX, where
two Binary32 floating-point additions are performed in parallel, in the same amount of time
as one Binary64 floating-point addition. 
The compiler may thus choose to execute $(a+b)+(c+d)$ in \texttt{Binary32},
where the two additions $a+b$ and $c+d$ are performed in parallel,
as it will execute faster than $((a+b)+c)+d$.
The compiler may also decide to use more accurate registers (64 bits or 80 bits) when such
vectorised devices are not available.
Thus, depending on the programming language, on the
compiler and its options and on the architecture, the intermediate results may vary,
as does the final result.
In some languages, it is possible to gain some a posteriori knowledge on the employed precision.
In C99, the value of \texttt{FLT\_EVAL\_METHOD} gives an indication, at run-time, on the
intermediate chosen format: indeterminate, double or long double.

\subsection{Problems for Interval Computations}
\label{SubSec.Prec.PbsIA}
The notion of precision in interval arithmetic could be regarded as the radius of the input arguments.
It is known that the overestimation of the result of a calculation is proportional to the 
radius of the input interval, with a proportionality constant which depends on the computed
expression.
More precisely \cite[Section 2.1]{Neumaier90:IA}, if $f$ is a Lipschitz-continuous function:
$D \subset \mathbb{R} ^n \rightarrow \mathbb{R}$, then
$$ q( \mathbf{f} (\mathbf{x}) , f(\mathbf{x})) = {\cal{O}} (r) \mbox{ if } \mathrm{rad}(\mathbf{x}) \leq r$$
where boldface letters denote interval quantities, 
$f(\mathbf{x})$ is the exact range of $f$ over the interval $\mathbf{x}$,
$\mathbf{f} (\mathbf{x})$ is an overestimation of $f(\mathbf{x})$ computed from an
expression for $f$ using interval arithmetic in a straightforward way,
and $q$ stands for the Hausdorff distance.
As the first interval here encloses the second, $q( \mathbf{f} (\mathbf{x}) , f(\mathbf{x})) $
is simply $\max ( \mathrm{inf}( f(\mathbf{x})) - \mathrm{inf} ( \mathbf{f} (\mathbf{x})) ,
                  \mathrm{sup} ( \mathbf{f} (\mathbf{x})) - \mathrm{sup}( f(\mathbf{x})) )$.
In this formula, $r$ can be considered as the precision.
It is possible to improve this result by using more elaborate approaches to interval arithmetic,
such as a Taylor expansion of order 1 of the function $f$.
Indeed, if $f$ is a continuously differentiable function then \cite[Section 2.3]{Neumaier90:IA},
$$ q( \mathbf{f} (\mathbf{x}) , f(\mathbf{x})) = {\cal{O}} (r^2) \mbox{ if } \mathrm{rad} (\mathbf{x}) \leq r$$
where $\mathbf{f} (\mathbf{x})$ is computed using a so-called centered form.
As $f$ is Lipschitz, it holds that the radius of $f (\mathbf{x})$
is also proportional to the radius of the input interval.
In other words, the accuracy on the result improves with the radius, or precision, of the input.
These results hold for an underlying exact arithmetic.
\\

This result in interval analysis can be seen as an equivalent to the rule of thumb
in numerical floating-point analysis \cite[Chapter 1]{Higham02}.
Considered with optimism, this means that it suffices to increase the computing
precision -- in floating-point arithmetic -- or the precision on the inputs -- in interval arithmetic --
to get more accurate results.
In interval arithmetic, this can be done through bisection of the inputs,
to get tight intervals as inputs and to reduce $r$ in the formula above.
The final result is then the union of the results computed for each subinterval.
A more pragmatic point of view is first that, even if the results are getting more and more
accurate as the precision increases, there is still no guarantee that the employed precision
will allow to reach a prescribed accuracy for the result.
Second, reproducibility is still out of reach, as developed in Section \ref{Sec.Order}.
\\

In what follows, we consider interval arithmetic implemented using floating-point arithmetic.
The aforementioned theorems, that bound $q( \mathbf{f} (\mathbf{x}) , f(\mathbf{x}))$,
do not account for the limited precision of the underlying floating-point
arithmetic. Indeed, computed interval results also  suffer from the problems of floating-point arithmetic,
namely
the possible loss of accuracy. Even if directed roundings make it possible to ensure the inclusion property,
there is no information about how much the exact interval is overestimated,
and no guarantee about the tightness of the computed interval.

\subsection{Recommendations}
\label{SubSec.Prec.Rec}
We advocate the use of the mid-rad representation on the one hand, and the use of iterative
refinement on the other hand.
The mid-rad representation $\langle m,r \rangle $ corresponds to the interval $[m-r, m+r] = \{ x \: : \: |m-x| \leq r \}$.
An advantage of the mid-rad representation is that thin intervals are represented more accurately
in floating-point arithmetic. For instance, let us consider $m$ a non-zero floating-point number
and $r = 1/2 \mathrm{ulp}(m)$, i.e. $r$ is a power of
2 that corresponds, roughly speaking, to half the last bit in the floating-point representation of $m$. 
(Let us recall \cite{Mu10:FP} that for $x \in \mathbb{R} \setminus \{0\}$,
if $x \in [2^e, 2^{e+1})$ then $\mathrm{ulp}(x) = 2^{e-p+1}$ in radix-2 floating-point arithmetic with $p$ bits used for the significand
and that $\mathrm{ulp}(-x) = \mathrm{ulp}(x)$ for negative $x$.)
In this example, $m$ and $r$ are floating-point numbers.
Then the floating-point representation by endpoints of the interval is $[ \mathrm{RD} (m-r), \mathrm{RU}(m+r)]$,
(where $\mathrm{RD}$ denotes the rounding mode towards $-\infty$ or rounding downwards
and $\mathrm{RU}$ denotes the rounding mode towards $+\infty$ or rounding upwards),
which is $[m-2r, m+2r]$:
in this example,
the width of the interval is doubled with the representation by endpoints.
The reader has to be aware that no representation, neither mid-rad nor by endpoints, always supersedes the other one.
An example where the mid-rad representation is superior to the representation by endpoints has just been given.
Conversely,
an unbounded interval can be represented by its endpoints, say $[1, + \infty)$, but the only enclosing
mid-rad representation $\langle m, + \infty \rangle $ with $m > 1$, corresponds to $\mathbb{R}$.

We also recommend the use of iterative refinement where applicable.
Indeed, even if the computations of the midpoint and the radius of the result suffer from the
aforementioned lack of accuracy, iterative refinement (usually a few iterations suffice) recovers
a more accurate result from this inaccurate one.

Let us illustrate this procedure on the example of square linear system solving $\mathbf{A} \mathbf{x} = \mathbf{b}$
(cf. \cite[Chapter 4]{Neumaier90:IA}, \cite[Chapter 7]{MoKeCl09:IA}, \cite{Rump05:IA} for an introduction).
Once an initial approximation $\mathbf{x^0}$ is computed, the residual $\mathbf{r} = \mathbf{b} - \mathbf{A} \mathbf{x^0}$
is computed \emph{using twice the current precision} (and here we rejoin the solutions in 
\cite{Bai12:repro,HeDi01:repro,RoRoAu11:repro,TaPaSaPa10-repro} already mentioned), as much cancellation occurs in this calculation.
Then, solve -- again approximately -- the linear system $\mathbf{A} \mathbf{e} = \mathbf{r}$ with the same matrix $\mathbf{A}$,
and re-use every pre-computations done on $\mathbf{A}$, typically a factorisation such as LU. Finally, correct the approximate solution: 
$\mathbf{x^1} \leftarrow \mathrm{mid} ( \mathbf{x^0} ) + \mathbf{e}$.
Under specific assumptions, but independently of the order of the operations, it is possible
to relegate the effects of the intermediate precision and of the condition number after the
significant bits in the destination format \cite{Nguyen11}.
In other words, the overestimation due to the floating-point arithmetic is minimal: only one ulp
(or very few ulps) on the precision of the interval result.

Another study \cite{Tisseur01} also takes into account the effect of floating-point arithmetic.
It suggests that the same approach applies to nonlinear system
solving and that the iterative refinement, called in this case Newton iteration,
again yields fully accurate results, i.e. up to 1 ulp of the exact result.

However, in both cases it is assumed that enough steps have been performed:
if there is a limit on the number of steps, then one run could converge but not the other one
and again numerical reproducibility would not be gained.

To conclude on the impact of the computing precision: it raises no problem for the validity of interval
computations, i.e. the inclusion property is satisfied, but it influences the accuracy of the result
and the execution time. It seems that the same could be said for the order of the operations,
which is the issue discussed next, but it will be seen that the validity of interval computations
can depend on it.

\section{Order of the Operations}
\label{Sec.Order}
\subsection{Problems for Interval Computations}
As already abundantly mentioned, a main explanation to the lack of reproducibility of
floating-point computations is the lack of associativity of floating-point operations
(addition, multiplication).

Interval arithmetic also suffers from a lack of algebraic properties.
In interval arithmetic, the square operation differs from the multiplication by the same
argument, because variables ($x$ and $y$ in the example) are decorrelated:
$$
\begin{array}{ll}
\multicolumn{2}{l}{[-1,2] ^2 = \{ x ^2, \: x \in [-1,2] \} =  [0,4]} \\
\neq & [-1,2] \cdot [-1,2]  \\
 & =  \{ x \cdot y, \: x \in [-1,2], \: y \in [-1,2] \} \\
 & =  [-2,4].
\end{array}
$$
This problem is often called \emph{variable dependency}.

In interval arithmetic, the multiplication is not distributive over the addition, again because of
the decorrelation of the variables.

Finally, interval arithmetic implemented using  floating-point arithmetic suffers from all these
algebraic features.
In any case, the computed result contains the exact result. However, it is not possible
to guarantee that one order produces tighter results than another one, as illustrated
by the following example.
In the Binary64 format, the smallest floating-point number larger than 1 is $1 + 2^{-52}$.
Let us consider three intervals with floating-point endpoints: $\mathbf{A}_1 = [- 2^{-53}, 2^{-52}]$,
$\mathbf{A}_2 = [-1, 2^{-52}]$ and $\mathbf{A}_3 = [1,2]$.
Using the representation by endpoints and the implementation
in floating-point arithmetic of $[\underline{a},\bar{a}] + [\underline{b},\bar{b}]$
as $[ \mathrm{RD}(\underline{a} + \underline{b}), \mathrm{RU} (\bar{a}+\bar{b})]$,
one gets for $(\mathbf{A}_1 + \mathbf{A}_2) + \mathbf{A}_3$:
$$
\begin{array}{lcl}
\mathbf{tmp}_1 & := & \mathbf{A}_1 + \mathbf{A}_2 \\
   & = & [ \mathrm{RD} (- 2^{-53} -1), \mathrm{RU}(2^{-52}+2^{-52})] \\
   & = & [-1-2^{-52}, 2^{-51}] \\
\multicolumn{3}{l}{\mbox{and finally}} \\
\mathbf{B}_1 & := & \mathbf{tmp}_1 + \mathbf{A}_3 \\
    & = & [ \mathrm{RD} (-1-2^{-52} + 1), \mathrm{RU}(2^{-51} + 2)] \\
    & = & [-2^{-52}, 2+2^{-51}].
\end{array}
$$
And one gets for $\mathbf{A}_1 + (\mathbf{A}_2 + \mathbf{A}_3)$:
$$
\begin{array}{lcl}
\mathbf{tmp}_2 & := & \mathbf{A}_2 + \mathbf{A}_3 \\
    & = & [ \mathrm{RD} (- 1 +1), \mathrm{RU}(2^{-52}+2)] \\
    & = & [0, 2+2^{-51}] \\
\multicolumn{3}{l}{\mbox{and eventually}} \\
\mathbf{B}_2 & := & \mathbf{A}_1 + \mathbf{tmp}_2 \\
\vspace*{-10pt} & & \\
    & = & [ \mathrm{RD} (-2^{-53} + 0), \mathrm{RU}(2^{-52} + 2+ 2^{-51} )] \\
    & = & [-2^{-53}, 2+2^{-50}].
\end{array}
$$
The exact result is $\mathbf{B} := [-2^{-53}, 2+2^{-51}]$ and this interval is representable using floating-point endpoints.
It can be observed that both $\mathbf{B}_1$ and $\mathbf{B}_2$ enclose $\mathbf{B}$: the inclusion property is satisfied.
Another observation is that $\mathbf{B}_1$ overestimates $\mathbf{B}$ to the left and $\mathbf{B}_2$ to the right.
Of course, one can construct examples where both endpoints are under- or over-estimated.
\\

Not only does the order
in which non-associative operations such as additions are performed matter, but other orders do as well.
Let us go back to the bisection process mentioned in Section \ref{SubSec.Prec.PbsIA}.
Indeterminism is present
in the bisection process, 
and it introduces more sources of irreproducibility.
Indeed, with bisection techniques, one interval is split into 2 
and one half (or the two halves) are stored in a list for later processing.
The order in which intervals are created and processed is usually not deterministic,
if the list is used by several threads. This can even influence the creation or bisection
of intervals later in the computation. Typically, in algorithms for global optimisation \cite{HaWa03:AO,ReDeMePl99:par},
the best enclosure found so far for the optimum can lead to the exploration or destruction
of intervals in the list.
As the order in which intervals are explored varies, this enclosure for the optimum
varies and thus the content of the list varies
-- and not only the order in which its content is processed.

\subsection{Recommendations}
\label{SubSec.Order.Rec}
As shown in the example of the sum of three intervals, even if the result depends
on the order of the operations, the inclusion property holds and it always encloses
the exact results.
Thus an optimistic view is that numerical reproducibility is irrelevant for
interval computations, as the result always satisfies the inclusion property.
An even more optimistic view could be that getting different results is good, since
the sought result lies in their intersection. Intersecting the various results would yield
even more accuracy. Indeed, with the example above, the exact result $\mathbf{B}$
is recovered by intersecting $\mathbf{B}_1$ and $\mathbf{B}_2$.
Unfortunately this is rarely the case, and usually no computation yields tighter results
than the other ones. Thus the order may matter.
\\

We have already mentioned the benefit of using the mid-rad representation in Section 
\ref{SubSec.Prec.Rec}.
Another, much more important, advantage is to be able to benefit from efforts made in developing
mathematical libraries.
As shown by Rump in his pioneering work \cite{Rump99},
in linear algebra in particular, it is possible to devise algorithms that use floating-point routines.
Typically these algorithms compute the midpoint of
the result using optimised numerical routines, and they compute afterwards the radius of the result
using again optimised numerical routines.
Let us illustrate this approach with the product of matrices with interval coefficients, given in \cite{Rump99}:
$\mathbf{A} = \langle A_m, A_r \rangle $ and $\mathbf{B} = \langle B_m, B_r \rangle $ are matrices with interval coefficients,
represented as matrices of midpoints and matrices of radii. The product $\mathbf{A} \cdot \mathbf{B}$ is
enclosed in $\mathbf{C} = \langle C_m, C_r \rangle $ where 
an interval enclosing $C_m$ is computed as $[\mathrm{RD} (A_m \cdot B_m), \mathrm{RU} (A_m \cdot B_m) ]$ using optimised BLAS3
routines for the product of floating-point matrices.
Then $C_r$ is computed using $\mathrm{RU} ( (|A_m| + A_r) \cdot B_r + A_r \cdot |B_m|)$,
again using optimised BLAS routines.
In \cite{Rump99}, the main benefit which is announced is the gain in execution time,
as these routines are really well optimised for a variety of architectures,
e.g. in Goto-BLAS or in ATLAS or in MKL, the library developed by Intel
for its architectures.
We also foresee the benefit of using reproducible libraries, once they are developed,
such as the CNR version of Intel MKL, and once their usage and performance
are optimised. 

To sum up, so  far only accuracy and speed of interval computations can be affected by the order
of operations, but not the validity of the results. As discussed now, the order may also impact
the validity of interval results, i.e. the inclusion property may not be satisfied.

\subsection{Concerns Regarding the Validity of Interval Results}
\label{SubSec.Order.Damper}
Apart from their lack of numerical reproducibility, there is another limit to the current usability of floating-point BLAS routines for interval
computations. This limit lies in the use of strong assumptions on the order in which operations are
performed. Let us again exemplify this issue on the product of matrices with interval coefficients.
The algorithm given above reduces interval operations to floating-point matrices operations.
However, it requires 4 calls to BLAS3 routines. An algorithm has recently been proposed \cite{Rump12}
that requires only 3 calls to floating-point matrix products.
An important assumption to ensure the inclusion property is that two of these floating-point matrix products
perform the basic arithmetic operations in the same order.
Namely, the algorithm relies on the following theorem \cite[Theorem 2.1]{Rump12}:
if $A$ and $B$ are two $n \times n$ matrices with floating-point coefficients,
if $C = \mathrm{RN}(A \cdot B)$ is computed in any order and 
if $\Gamma = \mathrm{RN}(|A| \cdot |B|)$ is computed \emph{in the same order},
then the error between $C$ and the exact product $A \cdot B$ satisfies
$$ | \mathrm{RN}(A \cdot B) - A \cdot B| \leq \mathrm{RN}((n+2)u \Gamma + \eta)$$
where $u$ is the unit roundoff and $\eta$ is the smallest normal positive floating-point number
($u = 2^{-52}$ and $\eta = 2^{-1022}$ in Binary64).
This assumption on the order is not guaranteed by any BLAS library we know,
and the discussion above tends to show that this assumption does not hold for multithreaded computations.
\\

Two solutions can be proposed to preserve the efficiency of this algorithm along with the inclusion
property.
A straightforward solution consists in using a larger bound for the error on the matrix product,
a bound which holds whatever the order for the computations of $C$ and $\Gamma$.
Following \cite[Chapter 2]{Higham02} and barring underflow, a bound of the following form
can be used instead \cite{Revol13}:
$$ \begin{array}{l}
 | \mathrm{RN} (A \cdot B) - A \cdot B| \leq  \\
\hspace*{-2mm} 
   \mathrm{RN} \left( \{ (1 \! + \! 3u)^n \! - \! 1 \} \! \cdot  \!
   \left( \{(1 \! + \! 2u)|A| \} \! \cdot \! \{ (1 \! + \! 2u) |B| \} \right) \right) .
\end{array} $$

A more demanding solution, implemented in \cite{ReTh13:IMM}, consists in implementing simultaneously the two
matrix products, in order to ensure the same order of operations.
Optimisations to get performances are done: block products to optimise the cache usage (Level 1)
and hand-made vectorisation of the code yield performances, even if the performances of well-known BLAS
are difficult to attain.
Furthermore, this hand-made implementation has a better control on the rounding modes
than the BLAS routines, and this is another major point in satisfying the inclusion property,
as discussed in the next Section.
The lack of respect of the rounding modes is another source of interval invalidity (to paraphrase
the term ``numerical reproducibility'') and it is specific to interval computations.

\section{Rounding Modes}
\label{Sec.Rndg}
\subsection{Problems for Interval Computations}
\label{SubSec.Rndg.Pbs}

Directed rounding modes of floating-point arithmetic are crucial for a correct implementation
of interval arithmetic, i.e. for an implementation that preserves the inclusion property.
This is a main difference between interval computations and floating-point computations,
that usually employ rounding-to-nearest only. It is also a issue not only for reproducibility,
but even for correctness of interval computations, especially for parallel implementations, as it will be
shown below.

For some floating-point computations that make extensive use of EFT (Error-Free Transforms),
such as compensated sums mentioned in Section \ref{Sec.Prev},
the already mentioned QD library \cite{HiLiBa:QD},
or for XBLAS, a library for BLAS with extended precision \cite{XBLAS02}, the only usable rounding mode is
rounding-to-nearest, otherwise these libraries fail to deliver meaningful results in higher precision.

For interval arithmetic, both endpoints and mid-rad representations require directed rounding modes.
However, many obstacles are encountered.
First, it may happen that the compiler is too eager to optimise a code to respect the required
rounding mode. For instance, let us compute the enclosure of $a/b$ where both $a$ and $b$ are
floating-point numbers, using the following piece of pseudo-code: 
\begin{verbatim}
set_rounding_mode (downwards)
left := a/b
set_rounding_mode (upwards)
right := a/b
\end{verbatim}
In order to spare a division
which is a relatively expensive operation, a compiler may 
assign \verb+right := left+, but then the result is not what expected.
This example is given in the gcc bug report \#34678 entitled
\emph{Optimization generates incorrect code with -frounding-math options}.
Even using the right set of compilation options does not solve the problem.
Usually these options, when doing properly the job, 
neutralise optimisations done by the compiler and the resulting
code may exhibit poor performances.

Second, interval computations may rely on floating-point routines, such as the BLAS routines 
in our example of the product of matrices with interval coefficients. For interval computations,
it is desirable that the
BLAS library respects the rounding mode set before the call to a routine.
(However C99 seems to exclude that external libraries called from C99 respect the rounding mode in use,
on the contrary C99 allows them to set and modify the rounding mode.)
This desirable behaviour is not documented for the libraries we know of, as
developed in \cite{LaMeMo12}. The opposite has been
observed by A. Neumaier and years later by F. Goualard in MatLab, as explained in his message
to \url{reliable_computing@interval.louisiana.edu} on 29 March 2012 quoted below.
\begin{quote}
In a 2008 mail from Pr. Neumaier to Cleve Moler forwarded to the IEEE
1788 standard mailing list
(\url{http://grouper.ieee.org/groups/1788/email/msg00204.html}), it is
stated that MATLAB resets the FPU control word after each external
call, which would preclude any use of fast directed rounding through
some C MEX file calling, say, the fesetround() function.

According to my own tests, the picture is more complicated than that.
With MATLAB R2010b under Linux/64 bits I am perfectly able to switch
the rounding direction in some MEX file for subsequent computation in
MATLAB. The catch is that the rounding direction may be reset by some
later events (calling some non-existent function, calling an M-file
function, ...). In addition, I cannot check the rounding direction
using fegetround() in a C MEX file because it always returns that the
rounding direction is set to nearest/even even when I can check by
other means that it is not the case in the current MATLAB environment.
\end{quote}

A more subtle argument against the naive use of mathematical libraries is based
on the monotony of the operations. Let us take again the example of a matrix product,
with nonnegative entries, such as in $\Gamma$ introduced in Section
\ref{SubSec.Order.Rec}.
If the rounding mode is set to rounding upwards for instance, then it suffices
to compute $\Gamma _{i,j}$ as $\mathrm{RU} (\sum_k |a_{i,k}| \cdot |b_{k,j}|)$
to get an overestimation of $\sum_k |a_{i,k}| \cdot |b_{k,j}|$. However, if $\Gamma$
is computed using Strassen's formulae, then terms of the form $x - z$ and $y + z $
are introduced. To get an upper bound on these terms, one needs an overestimation $\bar{x}$
of $x$, $\bar{y}$ of $y$ and $\bar{z}$ of $z$, but also an underestimation $\underline{z}$
of $z$. The overestimation of $x-z$ can thus be computed as $\mathrm{RU}(\bar{x} - \underline{z})$ and
the overestimation of $y+z$ as $\mathrm{RU}(\bar{y} + \bar{z})$.
In other words, one may need both over- and under-estimation of intermediate values,
i.e. one may need to compute many intermediate quantities twice. This would ruin
the fact that Strassen's method is a fast method. 
Anyway, if a library implementing Strassen's product is called with the rounding mode
set to $+ \infty$ and respects it, it simply performs all operations with this rounding
mode and there is no guarantee that the result overestimates the exact result.
To sum up, there is little evidence and little hope that mathematical libraries
return an overestimation of the result when the rounding mode is set to $+ \infty$.

Third, the programming language may or may not support changes of the rounding mode.
OpenCL does not support it, as explicitly stated in the documentation
(from the OpenCL Specification Version: 1.2, Document Revision: 15):
\emph{Round to nearest even is currently the only rounding mode required by the OpenCL
specification for single precision and double precision operations and is therefore
the default rounding mode. In addition, only static selection of rounding mode is supported.
Dynamically reconfiguring the rounding modes as specified by the IEEE 754 spec is unsupported.}
OpenMP does not support it either, as less explicitly stated in the documentation
(from OpenMP Application Program Interface, Version 4.0 - RC 2 - March 2013, Public Review Release Candidate 2):
\emph{This OpenMP API specification refers to ISO/IEC 1539-1:2004 as Fortran 2003. The
following features are not supported:
IEEE Arithmetic issues covered in Fortran 2003 Section 14 [\ldots]}
which are the issues related to rounding modes.

Fourth, the execution environment, i.e. the support for multithreaded execution,
usually does not document either how the rounding modes are handled.
For architectures or instructions sets which have the rounding mode specified in the code
for the operation, as CUDA for GPU, or as IA64 processors but without access from a high-level
programming language, rounding modes are handled smoothly, even if they are accessed
with more or less ease.
For other environments where the rounding modes are set via a global flag, it is not
clear how this flag is handled: it is expected that it is properly saved and restored
when a thread is preempted or migrated,
it is not documented whether concurrent threads on the same core ``share'' the rounding mode
or whether each of them can use its own rounding mode.
To quote \cite{LaMeMo12}:
\emph{How a rounding mode change is propagated
from one thread or node of a cluster to all others is unspecified in the C standard.
In MKL the rounding mode can be specified only in the VML (Vector Math Library) part
and any multi-threading or clustering behavior is not documented.}
\\

After this discussion, it may appear utopian to rely too much on directed
rounding modes to ensure that the inclusion property is satisfied.

\subsection{Recommendations}
Our main recommendation is to use bounds on roundoff errors rather than using directed rounding modes,
\emph{when it is not safe to do so}.
These bounds must be computable using floating-point arithmetic.
They must be computed using rounding-to-nearest which is most likely to be in use.
They are preferably independent of the rounding mode.
For instance, the bound given in Section \ref{SubSec.Order.Damper} can be made independent
of the rounding mode by replacing $u$, which is the roundoff unit for rounding-to-nearest,
by $2u$ which is an upper bound for any rounding mode \cite{Revol13}.
Another example of this approach to floating-point roundoff errors is the 
FI\_LIB library \cite{Filib}. It was also adopted to account for roundoff errors by the
COSY library \cite{RMB03}.
This approach is rather brute force, but it is robust to the change of rounding mode.
Furthermore, if the algorithm implemented in a numerical routine is known, it is possible
to use this numerical routine and to get an upper bound on its result.
For instance, the bound in Section \ref{SubSec.Order.Damper} on the error of a matrix product
can be obtained using any routine for the matrix product, as long as it does not use
a fast method such as Strassen's.

\section{Conclusion}
As developed in the preceding sections, obtaining numerical reproducibility
is a difficult task. Getting it can severely impede performances in terms of execution time.
It is thus worth checking whether numerical reproducibility is really needed,
or whether getting guarantees on the accuracy of the results suffices.
For instance, as stated in \cite{HeDi01:repro} about climate models:
\emph{It is known that there are multiple stable regions in phase
space $[\ldots]$ that the climate system could be attracted to. However, due
to the inherent chaotic nature of the numerical algorithms involved, it is feared
that slight changes during calculations could bring the system from one regime to
another.} In this example, qualitative information, such as the determination of the attractor for the
system under consideration, is more important than reproducibility.
This questioning goes further in a talk by Dongarra, similar to \cite{Dongarra13:repro-HPC},
where he advocated the quest for a guaranteed accuracy and the use of small computing precision,
such as single precision, rather than the quest for bit-to-bit reproducibility,
for speed and energy-consumption reasons.
\\

However, numerical reproducibility may be mandatory. In such a case,
our main recommendation to conclude this work is the following methodology.
\begin{itemize}
\item Firstly, develop interval algorithms that are based on well-esta\-blished
numerical bricks, so as
to benefit from their optimised implementation.
\item Second, convince developers and vendors of these bricks
to clearly specify their behaviour, especially what regards rounding modes.
\item If the second step fails,
replicate the work done for the optimisation of the
considered numerical bricks, to adapt them to the peculiarities and
requirements of the interval algorithm. 
A precursor to this recommendation is the recommendation in \cite{Marker:scalablemm}
that aims at easing such developments:
\emph{In order to achieve near-optimal performance, library developers must be given access to routines or kernels that provide computational- and utility-related functionality at a lower level than the customary BLAS interface.}
This would make possible to 
use the lower level bricks used in high-performance BLAS, e.g. computations at the level of the blocks and not of the entire matrix.
\item Get free from the rounding mode by bounding, roughly but robustly, errors with formulas independent of the rounding mode if needed.
\end{itemize}

Eventually, let us emphasise that the problem of numerical reproducibility
is different from the more general topic called \emph{Reproducible research in computer science},
which is for instance developed in a complete issue of the \textit{Computing in Science \& Engineering} magazine \cite{Cise09:repro}.
Reproducible research corresponds to the possibility to reproduce computational results by keeping track of the
code version, compiler version and options, input data\ldots used to produce the results that are
often only summed up, mostly as figures, in papers.
A possible solution is to adopt a ``versioning'' system not only for code files, but also
for compilation commands, data files, binary files\ldots
However, floating-point issues are usually not considered in the related publications;
they very probably should.

\section*{Acknowledgement}
This work is partly funded by the HPAC project of the French Agence Nationale de la Recherche (ANR 11 BS02 013).

\ifCLASSOPTIONcaptionsoff
  \newpage
\fi



%

%
\begin{IEEEbiography}[{\includegraphics[width=1in,height=1.25in,clip,keepaspectratio]{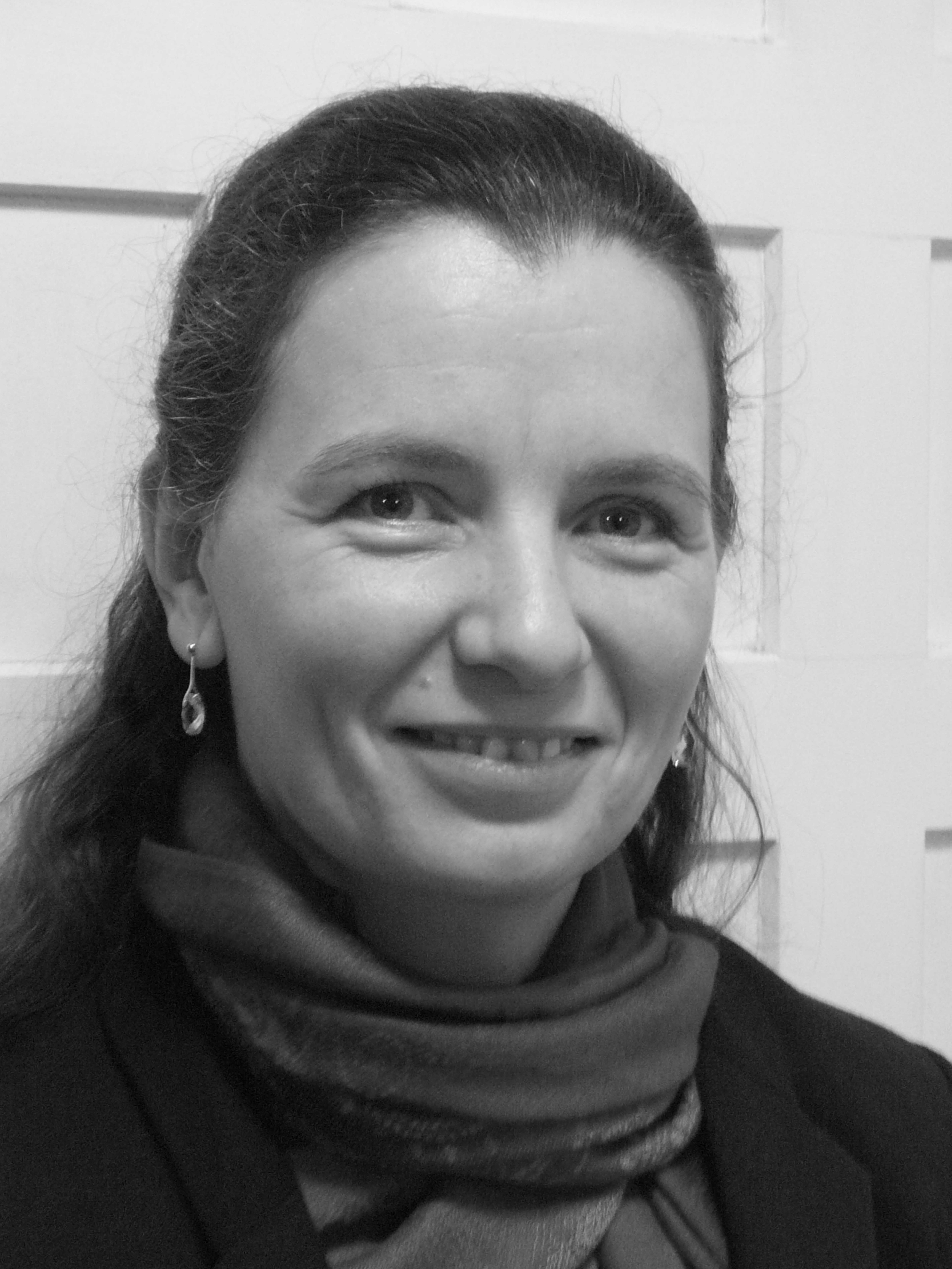}}]{Nathalie Revol}
received her PhD degree in applied mathematics from INPG (Grenoble, France) in 1994. She has been assistant
professor at the University of Sciences and Technologies of Lille, France, between 1995 and 2002.
Since then she is an INRIA researcher at the LIP laboratory, ENS de Lyon, France.
Since 2008 she co-chairs the IEEE P1788 working group for the standardization of interval
arithmetic with  R.~B. Kearfott.
Her research interests include interval arithmetic, floating-point arithmetic, and numerical quality.
\end{IEEEbiography}

\begin{IEEEbiography}[{\includegraphics[width=1in,height=1in,keepaspectratio]{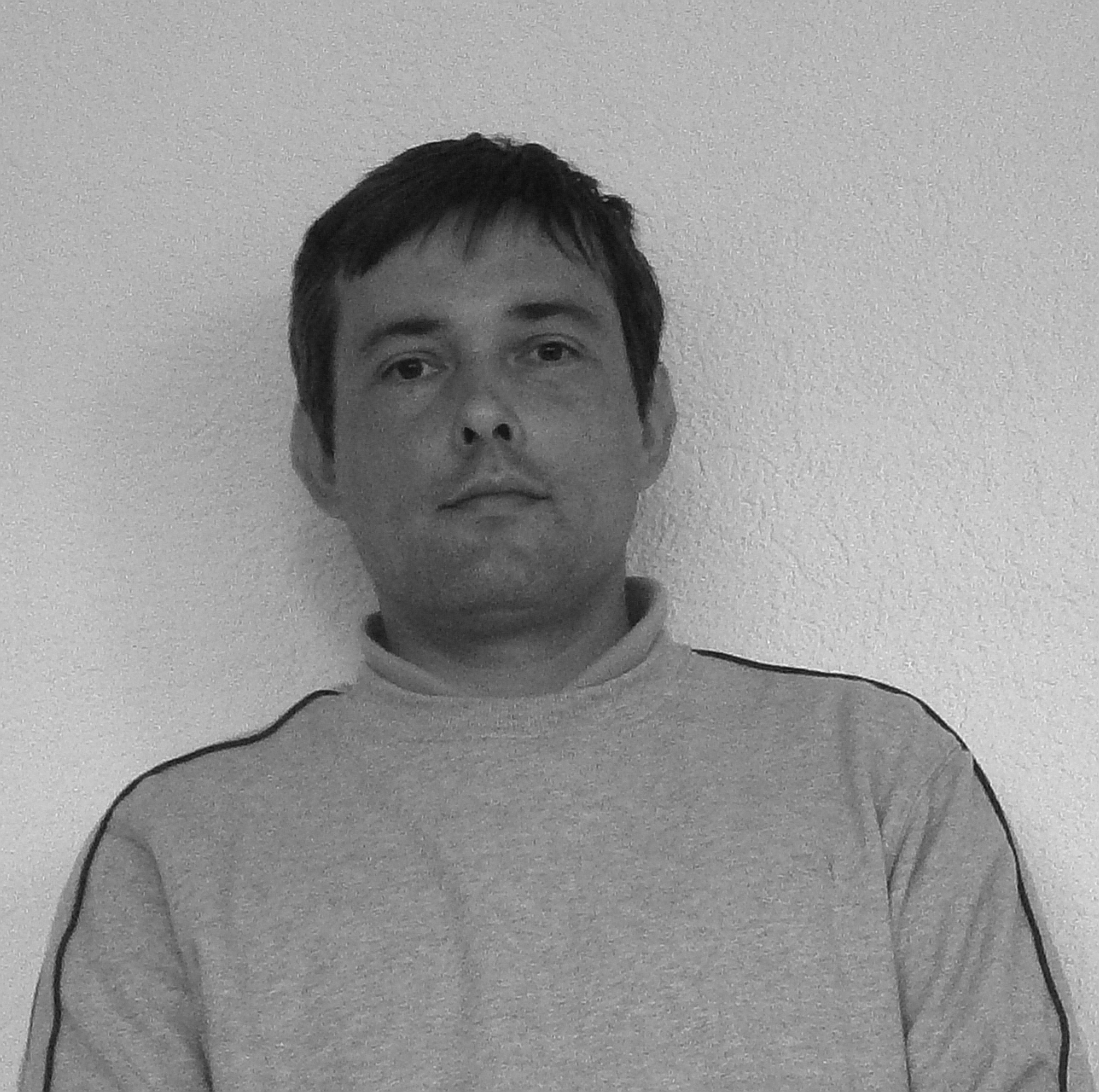}}]{Philippe Th\'eveny}
earned a MS degree in Computer Sciences from Universit\'{e} de Versailles
Saint-Quentin-en-Yvelines (France) in 2011. He is now a PhD candidate in
\'{E}cole Normale Sup\'{e}rieure de Lyon (France) under the supervision of
Nathalie Revol.
\end{IEEEbiography}





\end{document}